\documentclass[pre,onecolumn]{revtex4-2}
\usepackage{amsmath}  % needed for \tfrac, \bmatrix, etc.
\usepackage{amsfonts} % needed for bold Greek, Fraktur, and blackboard bold
\usepackage{graphicx} % needed for figures
\usepackage{subfigure}

\usepackage{tcolorbox}
\usepackage{alltt}

\begin{document}

\title{Root finding via local measurement}

\author{Jonathan Landy}
\affiliation{Stitch Fix, Inc.}
\email{jslandy@gmail.com}

\author{YongSeok Jho}
\affiliation{Gyeongsang National University}
\email{ysjho@gnu.ac.kr}

\begin{abstract}
We consider the problem of numerically identifying roots of a target function
-- under the constraint that we can only measure the derivatives of the
function at a given point, not the function itself. We describe and
characterize two methods for doing this:  (1) a local-inversion ``inching
process", where we use local measurements to repeatedly identify approximately
how far we need to move to drop the target function by the initial value over
$N$, an input parameter, and (2) an approximate Newton's method, where we
estimate the current function value at a given iteration via estimation of the
integral of the function's derivative, using $N$ samples.  When applicable,
both methods converge algebraically with $N$, with the power of convergence
increasing with the number of derivatives applied in the analysis.
\end{abstract}

\maketitle

\section{Introduction}

Here, we consider the general, constrained problem of finding a root of a
function $y$ that we cannot measure directly, but for which we can measure its
derivatives.  Specifically, we imagine that we are given a particular point
$(x_0, y_0)$ on the curve -- i.e., an initialization point, where the function
value is provided -- and we want to use this and available derivatives to find
a nearby root $x_*$ of $y$.  Standard methods for root finding -- e.g.,
Newton's method and the bisection method -- cannot be applied directly in this
case because each requires a function evaluation with each iterative search
step.  We suggest two methods for root finding that avoid this, both of which
involve making use of derivatives to estimate how much $y$ has changed as we
adjust $x$ to new points, hopefully closer to a root of $y$.  The first of
these, we call the local inversion method -- this proceeds by ``inching"
towards the root in $N$ steps, at each step adjusting $x$ by the amount needed
to approximately decrease $y$ by $y_0 / N$.  The second method, an approximate
Newton's method, makes use of approximate integrals of $y'$ to estimate the
function value at subsequent root estimates. 

In general, we expect the methods we discuss here to provide competitive
strategies whenever the target function $y$ -- whose root we seek -- is hard to
evaluate, relative to its derivatives.  This can of course occur in various
one-off situations, but there are also certain broad classes of problem where
this situation can arise.  For example, many physical systems are well-modeled
via differential equations that relate one or more derivatives of a target
function to a specified driving term -- Newton's laws of motion take this form.
Given an equation like this, we can often more readily solve for the
derivatives of the target than for the target itself -- whose evaluation may
require numerical integration. One simple example: Consider a flying, EV drone
whose battery drain rate is a given function of its instantaneous height $h$,
horizontal speed $\partial_t r$, and rate of climb $\partial_t h$.  If the
drone moves along a planned trajectory beginning with a full battery, root
finding techniques can be applied to the total charge function -- the integral
of the instantaneous rate of discharge -- to forecast when in the future the
battery will reach a particular level.  With the methods we describe here, this
can be done without having to repeatedly evaluate the precise charge at
candidate times.

Whereas standard methods for root finding -- e.g., Newton's method and the
bisection method --  generally converge exponentially quickly in the number of
function evaluations applied, the error in the methods we discuss here converge
to zero only algebraically.  However, the power at which we get convergence
increases with the number of derivatives supplied, and so we can often obtain
accurate estimates quite quickly.  Another virtue of the methods we discuss
here is that the simplest forms of the equations we consider are quite easy to
implement in code.  Given these points, these methods may be of general
interest for those who work with applied mathematics.

The remainder of our paper proceeds as follows:  In the next two sections, we
provide overviews of the two principal methods we discuss as well as a hybrid
strategy.  In the following section, we cover some example applications.  We
then conclude with a quick summary discussion.

Note:  We have written up the main strategies we discuss here in an open source
python package, \texttt{inchwormrf}.  This is available for download on pypi,
and on github at \texttt{github.com/efavdb/inchwormrf}.  The code can be
applied directly or used as demo code from which to base other implementations.

\section{Local inversion}

\subsection{Algorithm strategy}
\begin{figure}[h]\scalebox{.4}
{\includegraphics{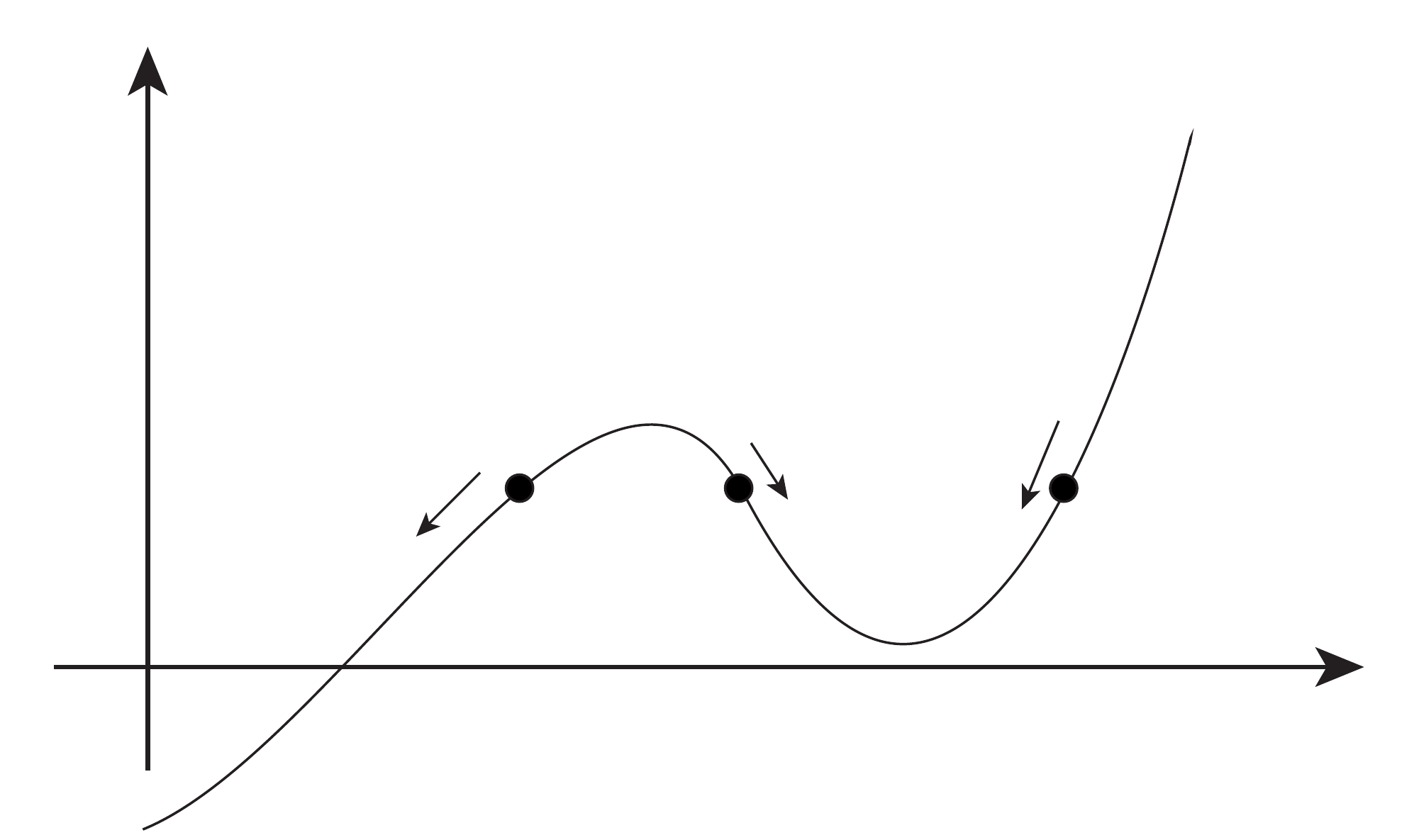}}
\caption{\label{fig:downhill} 
The local inversion method works its way ``downhill" and will find a root,
provided we can get to one without first encountering a zero derivative point
of $y$.  E.g., in this cartoon, a root will be found if we initialize at the
left-most black dot, but not if we initialize at either of the other two points
indicated. 
}
\end{figure}

The first strategy we will discuss is the local inversion method. We suppose we
measure $(x_0,y_0)$ and seek a root of $y$.  Our strategy will be to
iteratively work towards the root in $N$ steps, with $x_k$ chosen so that
$y(x_k) \approx y_0 (1 - k /N)$, so that $x_N$ is an approximate root.  At a
given step, then, the goal of the algorithm is to inch its way ``downhill" by
an amount $y_0 / N$, towards $y=0$.  This process will terminate near a root,
provided one does indeed sit downhill of the initial $(x_0, y_0)$ point --
i.e., provided one can get to a root without having to cross any zero
derivative points of $y$.  Fig.\ \ref{fig:downhill} illustrates the idea.

In order to accomplish the above, we need an approximate local inverse of $y$
that we can employ to identify how far to move in $x$ to drop $y$ by $y_0 /N$
at each step.  Here, we will focus on use of a Taylor series to approximate the
inverse, but other methods can also be applied.  To begin, we posit an
available and locally convergent Taylor series for $y$ at each $x_k$,
\begin{eqnarray}  \label{series}
dy \equiv\sum_{j=1}^m a_j d x^j + O \left( d x^{m+1} \right) 
\end{eqnarray}
We can invert this series to obtain an estimate for the change $dx$ needed to
have $y$ adjust by $dy$, as
\begin{eqnarray}\label{series_inversion}
dx =  \sum_{j=1}^m A_j dy^j + O \left( dy^{m+1} \right) 
\end{eqnarray}
where the coefficient $A_J$ is given by \cite{morsemethods},
\begin{equation}\label{inversion_coefs}
  A_n = \frac{1}{na_1^ n}\sum_{s_1, s_2, s_3,\cdots}(-1)^{s_1+s_2+s_3+\cdots}
\frac{n(n+1)\cdots(n-1+s_1+s_2+s_3+\cdots)}{s_1!s_2!s_3!\cdots}
\left(\frac{a_2}{a_1}\right)^{s_1}\left(\frac{a_3}{a_1}\right)^{s_2}\cdots.
\end{equation}
Here, the sum over $s$ values is restricted to partitions of $n-1$,
\begin{equation}
s_1+2s_2+3s_3+\cdots = n-1.
\end{equation}
The first few $A_i$ given by (\ref{inversion_coefs}) are,
\begin{eqnarray}
A_1 &=& \frac{1}{a_1} \\
A_2 &=& -\frac{1}{a_1^3} a_2 \\
A_3 &=& \frac{1}{a_1^5} \left(2 a_2^2 - a_1 a_3 \right).
\end{eqnarray}
In a first order approximation, we would use only the expression for $A_1$
above in (\ref{series_inversion}).  This would correspond to approximating the
function $y$ as linear about each $x_k$, using just the local slope to estimate
how far to move to get to $x_{k+1}$.  We can obtain a more refined estimate by
using the first two terms, which would correspond to approximating $y$ using a
local quadratic, and so on.  In practice, it is convenient to ``hard-code" the
expressions for the first few coefficients, but the number of terms $p(k)$
present in $A_k$ increases quickly with $k$, scaling as $\log p(k) \sim
\sqrt{n}$ \cite{andrews}.  For general code, then, it's useful to also develop
a method that can directly evaluate the sums (\ref{inversion_coefs}) for large
$n$.

This gives our first strategy: Begin at $(x_0, y_0)$, then write
\begin{eqnarray}
x_k = x_{k-1} + d x_{k-1}, \text{\ \ for $k \in \{1, \ldots, N \}$}
\end{eqnarray}
with $x_N$ being the final root estimate and $d x_{k-1}$ being the value
returned by (\ref{series_inversion}), in an expansion about $x_{k-1}$, with the
shift in $y$ set to $dy = -y_0 / N$.  Example code illustrating this approach
with only one derivative is shown in the boxed Algorithm 1.

\subsection{Asymptotic error analysis}

\begin{tcolorbox}[float, floatplacement=t]
\textbf{Algorithm 1: Local inversion (with m = 1 derivative only)}
\begin{verbatim}
def local_inversion_simple(x0, y0, dfunc, N):
    x_val = x0 * 1.0
    dy = y0 / float(N)
    for _ in range(N):
        x_val -= dy / dfunc(x_val)
    return x_val
\end{verbatim}
\end{tcolorbox}

Upon termination, the local inversion algorithm produces one root estimate,
$x_N \approx x_*$.  As noted above, one requirement for this to be a valid root
estimate is that there be a root downhill from $x_0$, not separated from it by
any zero derivative points of $y$.  A second condition is that the inverse
series (\ref{series_inversion}) provide reasonable estimates for the local
inverse at each $x_k$.  This should hold provided the inverse series
converges to the local inverse at each point, within some finite radius
sufficiently large to cover the distance between adjacent sample points.

We can identify a sufficient set of conditions for the above to hold as
follows:  If at the point $x$, $y$ is analytic, with (\ref{series}) converging
within a radius $R$ about $x$ and $\vert y(\tilde{x}) - y(x) \vert \leq M$ for
$\tilde{x}$ within this radius, and $y'(x) = s$ is non-zero at $x$, then the
inverse of $y$ will also be analytic at $x$, and its Taylor series
(\ref{series_inversion}) will converge within a radius of at least
\cite{redheffer}
\begin{eqnarray} \label{redhf}
r = \frac{1}{4} \frac{(s R)^2}{M}.
\end{eqnarray}
If the minimum of (\ref{redhf}) between $x_0$ and $x_*$ is some positive value
$\tilde{r} > 0$, it follows that each of inverse expansions applied in our
algorithm will converge between adjacent sample points, provided $N \geq \vert
y_0 \vert / \tilde{r}$.

Under the conditions noted above, we can readily characterize the $N$
dependence in the error in the root estimate $x_N$. Setting $dy =
-\frac{y_0}{N}$ and plugging into (\ref{series_inversion}), Taylor's theorem
tells us that the change in $x$ as $y$ moves by $dy$ will be given by
\begin{eqnarray}
d x_{k} = \sum_{j=1}^m A_j \left (-\frac{y_0}{N} \right )^j + O
\left(N^{-(m+1)} \right)
\end{eqnarray}
That is, when using the truncated series expansion to approximate how far we
must move to drop $y$ by $y_0 / N$, we'll be off by amount that scales like
$N^{-(m+1)}$.  After taking $N$ steps like this, the errors at each step will
accumulate to give a net error in the root estimate that scales like
\begin{eqnarray}\label{inverse_method_error}
x_N  \sim x_* +  O(N^{-m}).
\end{eqnarray}
That is, the local inversion root estimate's error converges to zero
algebraically with $N$, with a power that increases with $m$, the number of
derivatives used in the expansions.

\begin{figure}[h]\scalebox{.5}
{\includegraphics{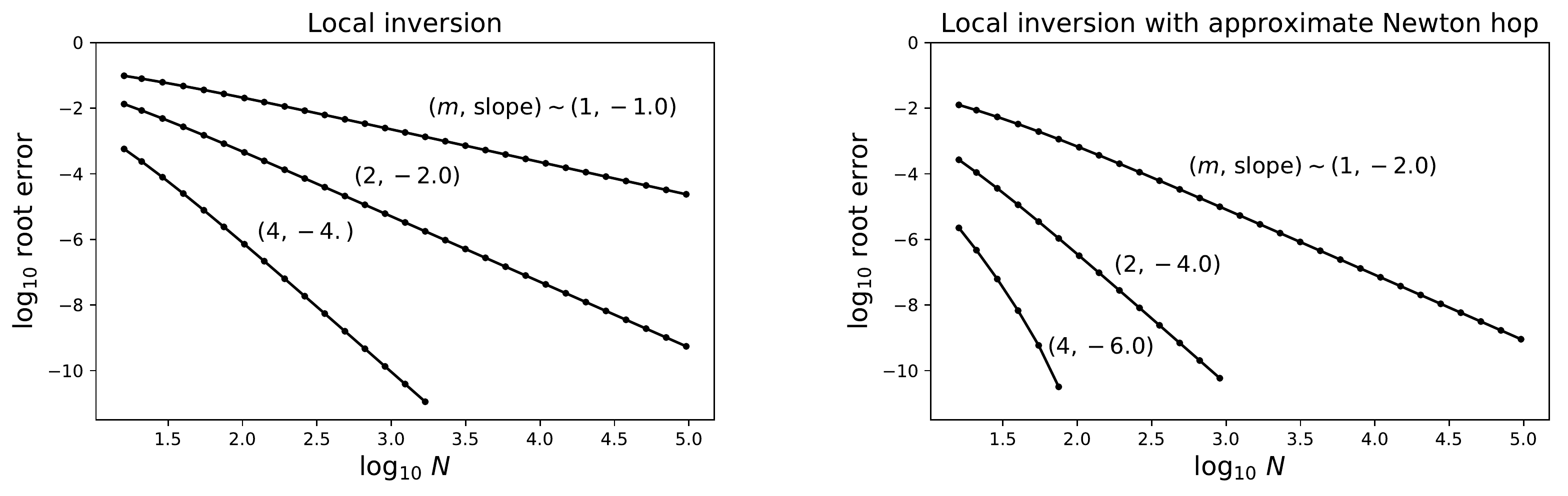}}
\caption{\label{fig:convergence} 
(left) Here, we plot the error in the root estimates obtained from the local
inversion method when applied to the function $y(x) = x^5 - 3$, the real root
of which is $3^{1/5} \approx 1.245731$.  We ran the algorithm initializing it
at the point $(x_0, y_0) = (2, 29)$, passing it $m=1$, $2$, or $4$ derivatives,
and obtained algebraic convergence to the correct root with errors scaling as
in (\ref{inverse_method_error}).  E.g., we get about four digits of accuracy if
we use just one derivative ($m=1$) and take $N=10^4$ steps, but six digits of
accuracy if we use four derivatives and only $N=10^2$ steps. (right)  Here, we
apply local inversion with a final, approximate Newton hop to the same problem
as above, again passing it $m=1$, $2$, or $4$, derivatives.  The errors
converge to zero more quickly now, in a manner consistent with
(\ref{Newton_error}).  Now, with four derivatives, we get around $12$ digits of
accuracy using $N=10^2$ steps. }
\end{figure}

The left plot of Fig.\ \ref{fig:convergence} illustrates this scaling in a
simple numerical experiment.  Here, we've applied the local inversion algorithm
to find the single real root of $y=x^5 - 3$. Plots of the error versus $N$ are
consistent with (\ref{inverse_method_error}) for each $m$.

\subsection{Formal expression for root}

If we take the limit of $N \to \infty$ of the $m=1$ version of the local
inverse algorithm, we obtain
\begin{eqnarray}\label{formal_expression}
x_* = x_0 - \int_0^{y_0} \frac{dy}{y'(x)},
\end{eqnarray}
a formal, implicit expression for the root we're estimating.  Note that this
depends only on the initial point $(x_0, y_0)$ and the derivative function,
$y'$.

Although we won't pursue this here in detail, we note that
(\ref{formal_expression}) provides a convenient starting point for identifying
many methods for estimating roots.  For example, we can recover Newton's method
by using a linear approximation to the integrand, expanding about $(x_0, y_0)$,
integrating this, then iterating.  Using a discrete approximation to the
integral, we can recover the local inversion algorithm above.  We can also
combine these ideas to obtain an iterated version of the local inversion
method.

A special class of functions for which (\ref{formal_expression}) can be quite
useful are those for which we can obtain an explicit expression for $y'$ in
terms of $y$.  In this case, (\ref{formal_expression}) becomes an explicit
integral which can be approximated in various ways, sometimes giving results
that converge much more quickly than does the general local inversion method.
For example, with $y = \cos(x)$ and $y' =- \sqrt{1-y^2}$,
(\ref{formal_expression}) returns a familiar integral for $\arccos$.   In a
numerical experiment, we applied Gaussian quadrature to obtain an estimate for
the integral from $y = 1/\sqrt{2}$ to $0$ -- giving estimates for $\pi / 4$.
The errors in these estimates converged exponentially to zero with the number
of samples $N$ used to estimate the integral -- beating the algebraic
convergence rate of local inversion.

\subsection{Behavior near a zero derivative}

When a zero derivative point separates $x_0$ and $x_*$, the local inversion
method will not return a valid root estimate.  However, we can still apply the
local inversion strategy if either $x_0$ or $x_*$ is itself a zero derivative
point.  The challenge in this case is that the inverse Taylor series
(\ref{series_inversion}) will not converge at the zero derivative point, and
some other method of approximating the inverse will be required near this
point.

To sketch how this can sometimes be done, we consider here the case where the
first derivative of $y$ is zero at $x_0$, but the second derivative is not.
This situation occurs in some physical applications.  E.g., it can arise when
seeking turning points of an oscillator, moving about a local minimum at $x_0$.
To estimate the change in $x$ needed for $y$ to move by $dy$ in this case, we
can apply the quadratic formula, which gives
\begin{eqnarray} \label{quadratic_formula}
d x \sim \frac{-a_1 + \sqrt{a_1^2 + 4 a_2 dy}}{2 a_2} + O(d x^3).
\end{eqnarray}
If we then run the local inversion algorithm with (\ref{quadratic_formula})
used in place of the inverse Taylor expansion, we can again obtain good root
estimates.  In this case, applying an error analysis like that we applied
above, one can see that the error here will be dominated by the region near the
zero derivative, giving a net root estimate error that scales like $N^{-3/2}$.
If there is no zero derivative point, we can still use
(\ref{quadratic_formula}) as an alternative to (\ref{series_inversion}).  In
this case, the root estimates that result again have errors going to zero as
$N^{-2}$ -- the rate of convergence that we get using series inversion with two
derivatives.

To improve on this approach using higher derivative information, one can simply
carry out a full asymptotic expansion of (\ref{series}) near the zero
derivative point.  This can be applied to obtain more refined estimates for the
$dx$ needed at each step in this region.  Away from the zero derivative point
this local asymptotic expansion will not converge, and one must switch over to
the inverse Taylor expansion (\ref{series_inversion}).  The resulting algorithm
is a little more complicated, but does allow one to obtain root estimates whose
errors go to zero more quickly with $N$.
\newline

This completes our overview of the local inversion root finding strategy.  We
now turn to our second approach, an approximate Newton hop strategy.

\section{approximate Newton hops}
\subsection{Basic strategy}
Here, we discuss our second method for root finding -- approximate Newton hops.
Recall that in the standard Newton's method, one iteratively improves upon a
root estimate by fitting a line to the function $y$ about the current position,
then moving to the root of that line.  That is, we take
\begin{eqnarray}\label{newtonsmethod}
x_{k+1} = x_k - \frac{y(x_k)}{y'(x_k)}
\end{eqnarray}
Iteratively applying (\ref{newtonsmethod}), the $x_k$ values often quickly
approach a root of the function $y$, with the error in the root estimate
converging exponentially quickly to zero.

We can't apply Newton's method directly under the condition that we consider
here, that $y$ is difficult to measure.  However, we assume that $y^{\prime}$
can be evaluated relatively quickly, and so we can apply Newton's method if we
replace the $y$ value in (\ref{newtonsmethod}) above by an approximate integral
of $y^{\prime}$, writing
\begin{eqnarray}\label{yprimeint}
y(x_k) = y(x_0) + \int_{x_0}^{x_k} y'(x) dx
\end{eqnarray} 
A simple strategy we can apply to evaluate the above is to sample $y^{\prime}$
at the $N$ equally-spaced points
\begin{eqnarray}
\tilde{x}_i = x_0 + \frac{i}{N-1} (x_{k} - x_0), \text{\ \ $i \in 0, \ldots, N
- 1$}.
\end{eqnarray}
Given samples of $y^{\prime}$ at these points, we can then apply the trapezoid
rule to estimate the integral in (\ref{yprimeint}).  This gives an estimate for
$y(x_k)$ that is accurate to $O(N^{-2})$.  Plugging this into
(\ref{newtonsmethod}), we can then obtain an approximate, Newton hop root
update estimate.  If we iterate, we then obtain an approximate Newton's method.
We find that this approach often brings the system to within the $O(N^{-2})$
``noise floor" of the root within $5$ to $10$ iterations.  If $y^{\prime} \not
= 0$ at the root, this then gives an estimate for the root that is also
accurate to $O(N^{-2})$.  A simple python implementation of this strategy is
given in the boxed Algorithm 2.

\begin{tcolorbox}[float, floatplacement=t]
\textbf{Algorithm 2: Approximate Newton hops (with m = 1 derivative only)}
\begin{verbatim}
def approximate_newton_simple(x0, y0, dfunc, N, iterations):
    x_val = x0 * 1.0
    
    for _ in range(iterations):
        # Step 1: Estimate y using approximate integral of y'
        y_val = y0 * 1.0
        dx = (x_val - x0) / float(N)
        for i in range(N + 1):
            x = x0 + i * dx
            y_val += dfunc(x) * dx
        y_val -= 0.5 * dfunc(x0) * dx
        y_val -= 0.5 * dfunc(x_val) * dx

        # Step 2: Newton hop
        x_val -= y_val / dfunc(x_val)
    
    return x_val
\end{verbatim}
\end{tcolorbox}

\subsection{Refinements}

We can obtain more quickly converging versions of the approximate Newton hop
strategy if we make use of higher order derivatives.  The Euler-Maclaurin
formula provides a convenient method for incorporating this information
\cite{apostol}.  We quote the formula below in a special limit, but first
discuss the convergence rates that result from its application:   If $m$
derivatives are supplied and we use $N$ points to estimate the integral at
right in (\ref{yprimeint}), this approach will give root estimates that
converge to the correct values as 
\begin{eqnarray}\label{Newton_error}
x_N \sim x_* + O(N^{-2 \lfloor \frac{m}{2} \rfloor - 2})
\end{eqnarray}
For $m=1$, the error term goes as $O(N^{-2})$, matching that of the trapezoid
rule we discussed just above.  For even $m \geq 2$, we get $O(N^{-m - 2})$.
That is, the error in this case goes to zero with two extra powers of $N$,
relative to that of the local inversion method, (\ref{inverse_method_error}).
This can give us a very significant speed up.

In our applications of the Euler-Maclaurin formula, we are often interested in a
slightly more general situation than that noted above, where the samples are
always equally spaced.  To that end, we'll posit now that we have samples of
$y^{\prime}$ at an ordered set of values $\tilde{x}_i$, with $i \in (0, \ldots,
N)$.  Direct application of the formula requires evenly spaced samples.  To
apply it in this case, then, we consider a change of variables, writing
\begin{eqnarray} \label{change_of_variables}
\int_{\tilde{x}_0}^{\tilde{x}_N} y^{\prime}(x) dx = \int_0^N
y^{\prime}(\tilde{x}(g)) \frac{d\tilde{x}}{dg} dg
\end{eqnarray}
Here, $\tilde{x}(g)$ is now some interpolating function that maps the domain
$[0, N]$ to $[\tilde{x}_0, \tilde{x}_N]$, satisfying
\begin{eqnarray}
\tilde{x}(g) = \tilde{x}_g, \text{ \  \ \ \ \ \ for $g \in \{0, 1, \ldots N\}$} 
\end{eqnarray}
With this change of variables, the integral at right in
(\ref{change_of_variables}) is now sampled at an evenly spaced set of points,
and we can use the Euler-Maclaurin formula to estimate its value.  If we plan
to use $m$ derivatives in our analysis, the interpolating function $\tilde{x}$
should also have at least $m$ continuous derivatives. In our open source
implementation, we have used a polynomial spline of degree $2 \lfloor m / 2
\rfloor + 1$ for $m> 1$ and a cubic spline for $m = 1$.

We now quote the Euler-Maclaurin formula, as applied to the right side of
(\ref{change_of_variables}).  This reads,
\begin{eqnarray} \label{em}
\int_{\tilde{x}_0}^{\tilde{x}_N} y^{\prime}(x) dx \sim \sum_{g=0}^N
y^{\prime}(\tilde{x}_g) \left .
\frac{d\tilde{x}}{dg} \right \vert_g - \frac{y^{\prime}(\tilde{x}_0) \left .
\frac{d\tilde{x}}{dg} \right \vert_{g=0} + y^{\prime}(x_N) \left .
\frac{d\tilde{x}}{dg} \right \vert_{g=N}}{2} + \sum_{k=1}^{\lfloor \frac{m}{2}
\rfloor} \frac{B_{2k}}{2k!} \left ( G^{(2k-1)}(\tilde{x}_N)  -
G^{(2k-1)}(\tilde{x}_0)\right)
\end{eqnarray}
Here $B_k$ is the $k$-th Bernoulli number and $G$ is the function of $g$ given
by
\begin{eqnarray}
G(g) = y^{\prime}(\tilde{x}(g))  \frac{d\tilde{x}}{dg}
\end{eqnarray}
Notice from (\ref{em}) that the higher order derivatives of $G$ are only needed
at the boundaries of integration, where $g=0, N$.   We can evaluate these
higher-order derivatives using the Faa di Bruno theorem \cite{roman}, which
generalizes the first derivative chain rule for derivatives of composite
functions. Each extra derivative with respect to $g$ gives another factor
scaling like $O(N^{-1})$, resulting in the scaling of the error quoted in
(\ref{Newton_error}) at fixed $m$.

With the formulas noted above, we can implement more refined versions of the
approximate Newton hop method: With each iteration, we apply (\ref{em}) to
approximate the current $y$ value, then plug this into (\ref{newtonsmethod}) to
obtain a new root estimate.  The samples used to estimate the integral of
$y^{\prime}$ can be equally spaced, or not.  Again, the former situation is
more convenient to code up, but the latter might be used to allow for more
efficient sampling:  E.g., one might re-use prior samples of $y^{\prime}$ in
each iteration, perhaps adding only a single new sample at the current root
estimate. Doing this, the samples will no longer each be evenly spaced, which
will force use of an interpolation function.  However, this can allow for a
significant speed up when $y^{\prime}$ evaluations are also expensive.

A hybrid local-inversion, final Newton hop approach can also be taken:
First carrying out the local inversion process with $m$ derivatives, will get
us to a root estimate that is within $O(N^{-m})$ of the true root.  If we save
the $y^{\prime}$ values observed throughout this process, we can then carry out
a single, final Newton hop using the estimate (\ref{em}) for the $y$ value at
termination of the local inversion process.  This will then bring us to a root
estimate accurate to $O(N^{-2 \lfloor \frac{m}{2} \rfloor - 2})$, as in
(\ref{Newton_error}).

In our open source package, we include a higher order Newton's method, but for
simplicity have only coded up the case where $N$ new samples of $y^{\prime}$
are evaluated with each iteration, always evenly spaced.  When $y^{\prime}$ is
easy to evaluate, this approach is much faster than that where a new
interpolation is required with each iteration.  However, to support the
situation where $y^{\prime}$ is also expensive to evaluate, we have also
implemented the hybrid local inversion-Newton method.  This requires only one
interpolation to be carried out in the final step, makes efficient use of
$y^{\prime}$ evaluations, and results in convergence rates consistent with our
analysis here.  The right plot of Fig.\ \ref{fig:convergence} illustrates this
rate of convergence in a simple application.

\section{Example applications}
\subsection{Functions with high curvature}
Here, we review the point that Newton's method will sometimes not converge for
functions that are highly curved, but that in cases like this the local
inversion method generally will.  Consider the function
\begin{eqnarray}
y = \vert x \vert^{\gamma}
\end{eqnarray}
where $\gamma > 0$ is a fixed parameter.  If we initialize Newton's method at
the point $(x_0, x_0 ^{\gamma})$, with $x_0
> 0$, a bit of algebra shows that Newton's method will move us to the point
\begin{eqnarray}
x_1 = x_0 \left( 1 - \frac{1}{\gamma} \right)
\end{eqnarray}
The distance of $x_1$ from the origin is a factor of $\vert 1- \gamma^{-1}
\vert$ times that of $x_0$. Repeated applications of Newton's method will
result in $\vert x_k \vert = x_0  \vert 1 - \gamma^{-1} \vert ^ k$. If $\gamma
> 1/2$, the method converges exponentially quickly to the root.  However, if
$\gamma < 1/2$, each iteration drives the estimate further away from the root
of $y$, due to ``overshooting".

High curvature is not an issue for the local inversion method:  By inching
along, it can quickly recalibrate to changes in slope, preventing overshooting.
Note, however, that applying a final, approximate Newton hop after the local
inversion process will no longer improve convergence in cases of high curvature
such as this.

\subsection{Smoothstep}
The smoothstep function \cite{bailey2009graphics}

\begin{eqnarray}
s_1(x) = 3 x^2 - 2 x^3, \ \ x \in (0, 1)
\end{eqnarray}
is an ``s-shaped" curve that is commonly used in computer graphics to generate
smooth transitions, both in shading applications and for specifying motions.
This function passes through $(0,0)$ and $(1,1)$ and has zero derivatives at
both of these points.  More generally, one can consider the degree-$(2n + 1)$
polynomial $s_n$ that passes through these points and has $n\geq 1$ zero
derivatives at both end points.  Simple inverse functions are not available for
$n>1$, so to determine where $s_n(x)$ takes on a specific value, one must often
resort to root finding.

\begin{figure}[h]\scalebox{.45}
{\includegraphics{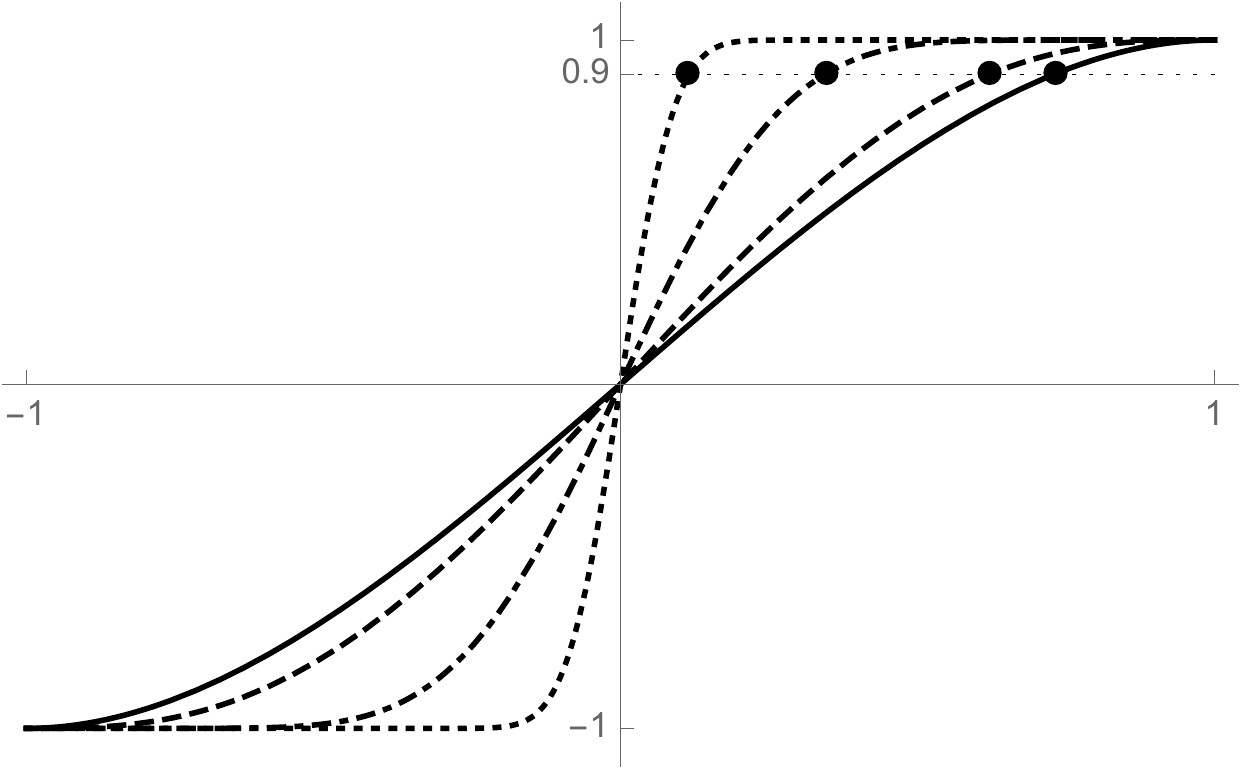}}
\caption{\label{fig:smoothstep} 
Here, we plot the $r_n(x)$ functions defined in (\ref{smoothstep_r}) for
$n = 1, 2, 10, 100$, and mark the points where each of these curves are
equal to $0.9$.  These values were found using the first and second derivatives
of $r$, initializing the approximate Newton method at the point $(0,-0.9)$.
}
\end{figure}

The following identity makes our present approach a tidy one for identifying
inverse values for general $n$: If we define,
\begin{eqnarray} \label{smoothstep_r}
r_n(x) = \frac{2 \Gamma \left(n+\frac{3}{2}\right)}{\sqrt{\pi } \Gamma (n+1)}
\int_0^x (1 - u^2)^n dx, \ \ x \in (-1, 1)
\end{eqnarray}
then
\begin{eqnarray}\label{smoothstep_rs}
s_n(x) = \frac{1}{2} r_n(2 x - 1)  - \frac{1}{2}.
\end{eqnarray}
Whereas evaluation of $r_n(x)$ requires us to evaluate the integral above,
the first derivative of $r_n$ is simply
\begin{eqnarray}
r_n^{(1)}(x) &=& \frac{2 \Gamma \left(n+\frac{3}{2}\right)}{\sqrt{\pi } \Gamma
(n+1)} (1 - x^2)^n 
\end{eqnarray}
Higher order derivatives of $r_n$ are easy to obtain from this last line as
well.  Given these expressions, we can easily invert $r_n$ for any value of $n$
using our methods.  Fig.\ \ref{fig:smoothstep} provides an illustration.

\subsection{Multiple dimensions}

Solutions to coupled sets of equations can also be treated using the local
measurement method.  We illustrate this here by example.  For simplicity, we'll
consider a cost function that can be written out by hand, but understand that
the method we describe will be most useful in cases where the derivatives of a
cost function are more easy to evaluate than the cost function itself.  The
function we'll consider is
\begin{eqnarray} \label{cost_function}
f(x_1, x_2) = x_1^4 + x_1^2 + 3 x_1 x_2 + x_2^2 + 7 x_1 + 9 x_2.
\end{eqnarray} 
If we wish to find a local extremum of this function, we can set the gradient
of the above to zero, which gives
\begin{eqnarray} \label{gradient} \nonumber
\nabla f &=& (4 x_1^3 + 2 x_1 + 3 x_2 + 7, 3 x_1 + 2 x_2 + 9) \\
&\equiv& (0, 0),
\end{eqnarray}
a coupled set of two equations for $x_1$ and $x_2$.

To apply the local inversion technique over $N$ steps, we need to expand the
gradient about a given position.  To first order, this gives
\begin{eqnarray} \label{2d_expansion}
d \nabla f = 
H \cdot (d x_1, d x_2)^T + O(d x^2),
\end{eqnarray}
where
\begin{eqnarray}
H = \left( \begin{array}{cc}
12 x_1^2 + 2 & 3  \\
3 & 2  \end{array} \right)
\end{eqnarray}
is the Hessian of (\ref{cost_function}), and
\begin{eqnarray} \label{target_grad_shift}
d \nabla f = - \frac{1}{N} \nabla f(x_{1,0}, x_{2,0}),
\end{eqnarray}
is the target gradient shift per step.  Inverting (\ref{2d_expansion}), we
obtain
\begin{eqnarray}\label{2d_inversion}
(d x_1, d x_2)^T = H^{-1} \cdot d \nabla f + O( d \nabla f^2).
\end{eqnarray}
There is a single real solution to (\ref{gradient}) at $(x_1, x_2) \approx (
1.35172698, -6.25699505)$.  Running the local inversion algorithm
(\ref{2d_inversion}) finds this solution with error decreasing like
$O(N^{-1})$.

To improve upon the above solution, we can either include more terms in the
expansion (\ref{2d_expansion}), or we can apply a final, approximate Newton
hop to refine the solution after the inching process.  To carry that out, we
need to consider how to estimate the final values of our gradient function.
To that end, we note that the integral we are concerned with is now a path
integral
\begin{eqnarray}\label{em_target_integral}
\nabla f = \nabla f(x_{1,0}, x_{2,0}) + \int_s H \cdot d \textbf{s}.
\end{eqnarray}
Taking a discrete approximation to this integral, and then applying a final
Newton hop, we obtain estimates with errors that converge to zero like
$O(N^{-2})$.

\section{Discussion}
Here, we have introduced two methods for evaluating roots of a function that
require only the ability to evaluate the derivative(s) of the function in
question, not the function itself.  This approach might be computationally
more convenient than standard methods whenever evaluation of the function
itself is relatively inconvenient for some reason.  The main cost of applying
these methods is that convergence is algebraic in the number of steps taken,
rather than exponential -- the typical convergence form of more familiar
methods for root finding, such as Newton's method and the bi-section method.
However, one virtue of these approaches is that higher order derivative
information is very easily incorporated, with each extra derivative provided
generally offering an increase to the convergence rate.  This property can
allow for the relatively slow rate of convergence concern to be mitigated.

The three main approaches that we have detailed here are as follows: (1)
the local inversion method, (2) the approximate Newton's method, and (3) the
hybrid local inversion-Newton method.  The first of these should be preferred
when looking at a function that is highly curved, since Newton's method is
subject to overshooting in this case.  The second will converge more quickly
for functions that are not too highly curved, and so is preferable in this
case.  Finally, the third method provides a convenient choice when working with
functions where Newton's method applies, but whose derivatives are also
somewhat costly to evaluate:  Our implementation here gives the convergence
rate of Newton's method and also minimizes the number of derivative calls
required.


\begin{thebibliography}{1}
\bibitem{bailey2009graphics}Bailey, M. and Cunningham, S.  Graphics shaders:
theory and practice. \textit{AK Peters/CRC Press} (2009)

\bibitem{morsemethods} Morse, P. M. and Feshbach, H.  Methods of Theoretical
Physics, Part 1 \textit{New York: McGraw-Hill} (1953)

\bibitem{andrews} Andrews, G. E. The Theory of Partitions \textit{Cambridge
University Press} (1976)

\bibitem{redheffer}Redheffer, R. M. Reversion of power series. Amer.  Math.
Monthly, 69(5):423–425 (1962)

\bibitem{roman} Roman, S. The formula of Faa di Bruno. Amer.  Math. Monthly
87.10:805-809 (1980)

\bibitem{apostol} Apostol, T. M. An Elementary View of Euler's Summation
Formula. Amer. Math. Monthly 106: 409-418 (1999).

\end{thebibliography}
\end{document}